\documentclass[10pt]{amsart}

\usepackage{amssymb,amsmath,amsthm}
\usepackage{graphicx,pgfplots}

\theoremstyle{plain}

\usepackage{graphicx,tikz}
\newtheorem{theorem}{Theorem}

\theoremstyle{definition}

\theoremstyle{remark}

\newcommand{\vol}{\operatorname{vol}}

\begin{document}

\title[]{Spectral Limitations of Quadrature Rules \\ and generalized spherical designs}

\author[]{Stefan Steinerberger}
\address{Department of Mathematics, Yale University}
\email{stefan.steinerberger@yale.edu}

\begin{abstract} We study manifolds $M$ equipped with a quadrature rule
$$ \int_{M}{\phi(x) dx} \simeq \sum_{i=1}^{n}{a_i \phi(x_i)}.$$
We show that $n-$point quadrature rules with nonnegative weights on a compact $d-$dimensional manifold  cannot integrate more than at most the first $c_{d}n + o(n)$ Laplacian eigenfunctions exactly. The constants $c_d$
are explicitly computed and $c_2 = 4$. The result is new even on $\mathbb{S}^2$ where it generalizes results on spherical designs.
\end{abstract}

\maketitle

\section{Introduction}
\subsection{Exact integration.} We study quadrature formulas
$$ \int_{M}{f(x) dx} \simeq \sum_{i=1}^{n}{a_i f(x_i)},$$
where $(M,g)$ is a compact Riemannian manifold without boundary, $(x_i)_{i=1}^{n} \subset M$ is a set of points and $(a_i)_{i=1}^{n} \in \mathbb{R}^{n}_{\geq 0}$ is a set of nonnegative weights. 
The fundamental question is how to choose the points and weights: a classical approach is to fix a set of functions and choose
points and weights so that these functions are integrated exactly; the canonical choice for functions on manifolds are the eigenfunctions of the Laplace operator
$$ -\Delta \phi_k = \lambda_k \phi_k.$$
 Any error in the quadrature formula is then a result of the presence high-frequency oscillation, which in a certain sense, cannot be avoided.
On the sphere, this idea goes back at least to a 1962 paper of Sobolev \cite{sobolev} (although he arrived there by different reasoning); the idea was then taken up by Lebedev \cite{leb1, leb5}
and is now sometimes called Lebedev quadrature. One would expect that $n$ points on $\mathbb{S}^2$, each having two coordinates and one
weight, should be able to integrate roughly the first $3n$ spherical harmonics exactly. This seems be a very good predictor (see Ahrens \& Beylkin \cite{ahrens}),
however, already McLaren \cite{mclaren} noted the existence of a set of $n=72$ points that integrates all polynomials up to the 14th degree exactly (corresponding
to $225 > 216 = 3\cdot 72$ functions). A natural question is whether, for $n$ sufficiently large, it is impossible to integrate more than $(3+\varepsilon)n$ spherical
harmonics exactly. The special case of all weights $a_i$ being identical and $M = \mathbb{S}^{d-1}$ is related to $t-$designs which have received a great deal
of attention resulting in hundreds of contributions, we refer to Delsarte, Goethals \& Seidel \cite{delsarte}, Seymour \& Zaslavsky \cite{sey}, Yudin \cite{yudin}, Bondarenko, Radchenko \& Viazovska \cite{bond, bond2}, the book of Conway \& Sloane \cite{conway} and references therein.
There also exists broad overlap with an entirely different set of problems; a 1904 question of J. J. Thomson \cite{thomson} is how
to distribute $n$ points on the sphere so as to minimize the energy
$$ \sum_{i,j=1 \atop i \neq j}^n{\frac{1}{\|x_i - x_j\|}}.$$
This line of question has been extended to general kernels $k(x_i, x_j)$ and studied in a very large number of settings, we refer to a recent survey of Brauchart \& Grabner \cite{brau}. One expects minimizers of these functionals to be spread very evenly over the manifold and for this reason minimizing configurations of various energy functionals are often used as sampling points (see \cite{brau}).

\subsection{Main results.} The purpose of our paper is to describe a very general inequality that connects these two different approaches (enforcing exactness on a finite-dimensional subspace
versus the notion of minimizing an energy functionals): we determine a special energy functional and prove that it is large whenever the quadrature formula is integrating a large number of
Laplacian eigenvalues exactly. As for notation, the Laplacian eigenfunctions $ -\Delta \phi_k = \lambda_k \phi_k$ are
assumed to be normalized in $L^2(M)$ and indexed so that $\phi_0$ is constant and $\phi_1$ is the first nontrivial eigenfunction. We will use $e^{t\Delta} f$ to denote the heat semigroup
$$ (\partial_t - \Delta) e^{t\Delta} f = 0$$
with $e^{0\Delta}f = f$. The Dirac measure in a point $x$ is denoted by $\delta_x$. The energy functional will be written in terms as pairwise interaction of Dirac measures
mollified via the heat
kernel; a useful heuristic is given by the short-time asymptotics
$$ \left[e^{t\Delta}\delta_x\right](y) \sim \frac{1}{(4\pi t)^{d/2}}\exp\left(-\frac{d(x,y)^2}{4t}\right).$$

\begin{center}
\begin{figure}[h!]
\begin{tikzpicture}[scale=0.7]
    \draw[thick, smooth cycle,tension=.7] plot coordinates{(-1,0) (0.5,2) (2,2) (4,3) (4.5,0)};
    \coordinate (A) at (1,1);
    \draw[thick] (A) arc(140:40:1) (A) arc(-140:-20:1) (A) arc(-140:-160:1);
\filldraw (4,1) ellipse (0.05cm and 0.05cm);
\draw[rotate around={45:(4,1)}, thick] (4,1) ellipse (1*0.28cm and 1*0.18cm);
\draw[rotate around={45:(4,1)}] (4,1) ellipse (2*0.28cm and 2*0.18cm);
\draw[rotate around={45:(4,1)}, dashed] (4,1) ellipse (3*0.28cm and 3*0.18cm);
\filldraw (0,0) circle (0.01cm);
\end{tikzpicture}
\caption{$\left[ e^{t\Delta} \delta_x\right](y)$ behaves roughly like a Gaussian at scale $\sim \sqrt{t}$ centered at $x$ whenever $t$ is small.}
\end{figure}
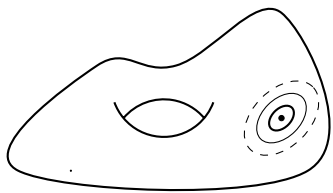
\end{center}

\begin{theorem} Let $(M,g)$ be a compact Riemannian manifold without boundary, let $(x_i)_{i=1}^{n} \subset M$ and $(a_i)_{i=1}^{n} \in \mathbb{R}^n_{>0}$ be given and
so that for all $0 \leq j \leq k$ the quadrature formula is exact for $\phi_j$.
Then, for some constant $c_{} > 0$ depending only on $M$, and all $t>0$
$$ \lambda_k \leq c  \frac{ \sum_{i,j=1}^{n}{a_i a_j\left(\frac{1}{t} + \frac{ d(x_i, x_j)^2}{t^2}\right) \left\langle e^{t\Delta} \delta_{x_i}, e^{t \Delta} \delta_{x_j} \right\rangle } }{  \left( \sum_{i,j=1}^{n}{a_i a_j \left\langle e^{t\Delta} \delta_{x_i}, e^{t \Delta} \delta_{x_j} \right\rangle }\right) - \vol(M)}$$\\
\end{theorem}

The inequality is close to sharp, especially for large values of $t$.
The right-hand side of the inequality can be thought of as a pairwise interaction functional. 
Since the bulk of the interactions is local and happens at scale $\sim \sqrt{t}$, we may, for $t$ small, replace the heat kernel by the short-time asymptotics to obtain an approximating energy that only depends on pairwise distances
$$ \mbox{ energy} = \frac{ \sum_{i,j=1}^{n}{\frac{a_i a_j}{(8 \pi t)^{d/2} }    \left(\frac{1}{t} + \frac{ d(x_i, x_j)^2}{t^2}\right)\exp\left(-\frac{d(x_i, x_j)^2}{8t}\right) }}{  \left( \sum_{i,j=1}^{n}{ \frac{a_i a_j}{(8 \pi t)^{d/2} }
\exp\left(-\frac{d(x_i, x_j)^2}{8t} \right)} \right) - \vol(M)}.$$
Our result guarantees that a good quadrature will have the quantity on the right-hand side be large. It is thus conceivable that actively
maximizing the functional may be a method to produce good quadrature points; in many instance it may be simpler to work with the even further
$$\mbox{simplified energy} = \sum_{i,j=1}^{n}{\frac{a_i a_j}{(8\pi t)^{d/2}} \exp\left( -\frac{d(x_i, x_j)^2}{8t}\right) }.$$

\begin{center}
\begin{figure}[h!]
\begin{tikzpicture}[scale=4]
\draw [thick] (0,0) -- (1,0) -- (1,1) -- (0,1) -- (0,0);
\filldraw (0.86344, 0.935215) circle (0.01cm);
\filldraw (0.761266, 0.123258) circle (0.01cm);
 \filldraw (0.312914, 0.914013) circle (0.01cm);
 \filldraw (0.282995, 0.136372) circle (0.01cm);
 \filldraw (0.520495, 0.702722) circle (0.01cm);
 \filldraw (0.889789, 0.411789) circle (0.01cm);
 \filldraw (0.38007, 0.420294) circle (0.01cm);
 \filldraw (0.828529, 0.685654) circle (0.01cm);
 \filldraw (0.0166807, 0.680519) circle (0.01cm);
 \filldraw (0.817869, 0.308049) circle (0.01cm);
 \filldraw (0.957418, 0.343145) circle (0.01cm);
 \filldraw (0.54677, 0.184955) circle (0.01cm);
 \filldraw (0.178159, 0.415532) circle (0.01cm);
 \filldraw (0.693605, 0.0526631) circle (0.01cm);
 \filldraw (0.715936, 0.639526) circle (0.01cm);
 \filldraw (0.441267, 0.982378) circle (0.01cm);
 \filldraw (0.222136, 0.92647) circle (0.01cm);
 \filldraw (0.468352, 0.596433) circle (0.01cm);
 \filldraw (0.580453, 0.469036) circle (0.01cm);
 \filldraw (0.788105, 0.875356) circle (0.01cm);
 \filldraw (0.607445, 0.781374) circle (0.01cm);
 \filldraw (0.025742, 0.575613) circle (0.01cm);
 \filldraw (0.607019, 0.373251) circle (0.01cm);
 \filldraw (0.275914, 0.412401) circle (0.01cm);
\filldraw (0.572954, 0.904129) circle (0.01cm);
 \filldraw (0.669615, 0.168117) circle (0.01cm);
 \filldraw (0.907688,  0.263678) circle (0.01cm);
 \filldraw (0.247459, 0.306925) circle (0.01cm);
 \filldraw (0.683555, 0.474396) circle (0.01cm); 
\filldraw (0.697176, 0.731277) circle (0.01cm);
 \filldraw (0.0779352, 0.0604842) circle (0.01cm);
 \filldraw (0.880141, 0.815375) circle (0.01cm);
 \filldraw (0.339601, 0.317686) circle (0.01cm); 
\filldraw (0.0129042, 0.171138) circle (0.01cm);
 \filldraw (0.68092, 0.842863) circle (0.01cm); 
\filldraw (0.522623, 0.828822) circle (0.01cm);
\filldraw (0.122374, 0.954634) circle (0.01cm);
 \filldraw (0.871111, 0.521851) circle (0.01cm);
 \filldraw (0.259988, 0.695398) circle (0.01cm);
 \filldraw (0.62305, 0.261789) circle (0.01cm);
 \filldraw (0.475526, 0.103445) circle (0.01cm);
 \filldraw (0.920645, 0.619608) circle (0.01cm);
 \filldraw (0.384579, 0.10938) circle (0.01cm); 
\filldraw (0.000487339, 0.793214) circle (0.01cm);
 \filldraw (0.06903, 0.874934) circle (0.01cm);
 \filldraw (0.108716, 0.628824) circle (0.01cm);
 \filldraw (0.354928, 0.00653162) circle (0.01cm); 
\filldraw (0.44633, 0.893699) circle (0.01cm);
 \filldraw (0.269775, 0.039552) circle (0.01cm);
 \filldraw (0.495805, 0.39909) circle (0.01cm);
 \filldraw (0.0867259,0.35875) circle (0.01cm); 
\filldraw (0.916754, 0.0126726) circle (0.01cm);
 \filldraw (0.594041, 0.0992317) circle (0.01cm);
 \filldraw (0.012737, 0.424625) circle (0.01cm);
 \filldraw (0.18215, 0.0362856) circle (0.01cm); 
\filldraw (0.921154, 0.724508) circle (0.01cm);
 \filldraw (0.154361, 0.130826) circle (0.01cm);
 \filldraw (0.986316, 0.0813935) circle (0.01cm); 
\filldraw (0.85324, 0.0845723) circle (0.01cm);
 \filldraw (0.448039,0.212008) circle (0.01cm);
 \filldraw (0.524963, 0.0132488) circle (0.01cm);
 \filldraw (0.0133802, 0.968421) circle (0.01cm);
 \filldraw (0.41018, 0.66293) circle (0.01cm);
 \filldraw (0.168091, 0.836753) circle (0.01cm);
 \filldraw (0.960149, 0.506272) circle (0.01cm); 
\filldraw (0.154239, 0.299023) circle (0.01cm); 
\filldraw (0.772125, 0.49642) circle (0.01cm);
 \filldraw (0.222319, 0.602971) circle (0.01cm);
 \filldraw (0.531493, 0.295142) circle (0.01cm);
 \filldraw (0.265543, 0.504527) circle (0.01cm);
 \filldraw (0.909212, 0.167291) circle (0.01cm);
 \filldraw (0.72026, 0.264402) circle (0.01cm);
 \filldraw (0.228343, 0.208178) circle (0.01cm);
 \filldraw (0.796125, 0.399643) circle (0.01cm);
 \filldraw (0.658114, 0.568048) circle (0.01cm);
 \filldraw (0.601584, 0.66728) circle (0.01cm); 
\filldraw (0.783496, 0.00141272) circle (0.01cm);
 \filldraw (0.433911, 0.317973) circle (0.01cm);
 \filldraw (0.159228, 0.53165) circle (0.01cm); 
\filldraw (0.328435, 0.605131) circle (0.01cm); 
\filldraw (0.349086,  0.217074) circle (0.01cm); 
\filldraw (0.107937, 0.207642) circle (0.01cm);
 \filldraw (0.617916, 0.986284) circle (0.01cm);
 \filldraw (0.362148, 0.838149) circle (0.01cm);
 \filldraw (0.445535, 0.768103) circle (0.01cm);
 \filldraw (0.700766, 0.366559) circle (0.01cm);
 \filldraw (0.806463, 0.204902) circle (0.01cm);
 \filldraw (0.0918888, 0.761954) circle (0.01cm);
 \filldraw (0.0242553, 0.269427) circle (0.01cm);
 \filldraw (0.468364, 0.491061) circle (0.01cm);
 \filldraw (0.783985, 0.773086) circle (0.01cm);
 \filldraw (0.253241, 0.805743) circle (0.01cm);
 \filldraw (0.37167, 0.521583) circle (0.01cm);
 \filldraw (0.16651, 0.71135) circle (0.01cm);
 \filldraw (0.798919, 0.593356) circle (0.01cm);
 \filldraw (0.701959, 0.933628) circle (0.01cm);
 \filldraw (0.342613, 0.737541) circle (0.01cm);
 \filldraw (0.0892862, 0.48232) circle (0.01cm);
 \filldraw (0.556831, 0.568482) circle (0.01cm); 
 \filldraw (0.954616, 0.885814) circle (0.01cm);
\end{tikzpicture}
\caption{A local minimizer of the simplified energy functional for $n=100$ points on $\mathbb{T}^2$ (found by starting with randomly chosen points and using gradient descent).}
\end{figure}
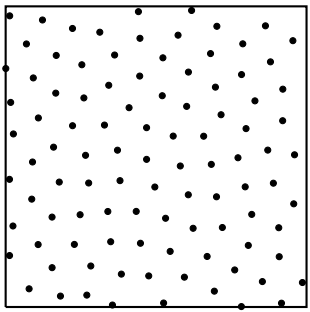
\end{center}

 One particularly important consequence is that $n$ points cannot integrate more than $\sim c_d n$ Laplacian eigenfunctions exactly. 

\begin{theorem} Let $M$ be a $d-$dimensional compact Riemannian manifold without boundary, let $(x_i)_{i=1}^n \subset M$ and $(a_i)_{i=1}^{n} \subset \mathbb{R}_{\geq 0}^n$. Then the quadrature rule
$$ \int_{M}{\phi(x) dx} \simeq \sum_{i=1}^{n}{a_i \phi(x_i)}$$
can integrate at most the first $c_d n + o(n)$ Laplacian eigenfunctions exactly, where $c_d$ depends only on the dimension. Moreover, $c_1 \sim 2.07, c_2=4$ and more generally
$$ c_d = \frac{\left( \frac{d}{2} + 1\right)^{\frac{d}{2} + 1}}{\Gamma\left(\frac{d}{2} + 1\right)} .$$
\end{theorem}
It is easy to see that $n$ equispaced points on the one-dimensional torus $\mathbb{T}$ can integrate the first $\sim 2n$  eigenfunctions exactly. 
Experimentally, there exist quadrature rules on $n$ points on $\mathbb{S}^2$ that integrate the first $\sim 3n$ spherical harmonics exactly \cite{ahrens}. This would
indicate that Corollary 2 is fairly sharp in one and two dimensions. The constants grow quickly ($c_3 \sim 7.43, c_4 = 13.5, c_5 \sim 24.13$). It would be interesting to understand whether the optimal value of $c_d$ grows linearly in the dimension and whether $c_d = d+1$. \\

We emphasize that the approach does not require the limit $n \rightarrow \infty$ to yield quantitative results: the method is equally applicable for finite $n \in \mathbb{N}$
and, especially when coupled with computational tools and some numerical analysis, should be able to provide fairly sharp computer-aided results for fixed values of $n$.
We also remark that the case of manifolds with boundary can be attacked by the same technique: there are some minor differences in terms of how the heat-flow is affected by the boundary
but these do not affect the method at large.  A variant of the method can even be applied to the case of finite graphs \cite{steingraph}.

\section{Proof of Theorem 1}
\begin{proof} It is an easy fact that for any function of the type
$$ g(x) = \sum_{\ell =k}^{\infty}{ \left\langle g, \phi_{\ell}\right\rangle \phi_{\ell}},$$
we have
\begin{align*}
 \left\| \nabla g\right\|_{L^2}^2 &= \left\langle \nabla \sum_{\ell =k}^{\infty}{ \left\langle g, \phi_{\ell}\right\rangle  \phi_{\ell}},  \nabla \sum_{\ell =k}^{\infty}{ \left\langle g, \phi_{\ell}\right\rangle  \phi_{\ell}} \right\rangle = \left\langle -\Delta \sum_{\ell =k}^{\infty}{ \left\langle g, \phi_{\ell}\right\rangle  \phi_{\ell}},  \sum_{\ell =k}^{\infty}{ \left\langle g, \phi_{\ell}\right\rangle  \phi_{\ell}} \right\rangle\\
&= \left\langle  \sum_{\ell =k}^{\infty}{ \lambda_{\ell} \left\langle g, \phi_{\ell}\right\rangle  \phi_{\ell}},  \sum_{\ell =k}^{\infty}{ \left\langle g, \phi_{\ell}\right\rangle  \phi_{\ell}} \right\rangle= \sum_{\ell=k}^{\infty}{\lambda_{\ell} \left|\left\langle g, \phi_{\ell}\right\rangle\right|^2 }\\
&\geq  \sum_{\ell =k}^{\infty}{\lambda_{k} \left|\left\langle g, \phi_{\ell}\right\rangle\right|^2 } = \lambda_k \left\|g\right\|_{L^2}^2.
\end{align*}
Suppose now that $(x_i)_{i=1}^{n} \subset M$ and $(a_i)_{i=1}^{n} \subset \mathbb{R}_{\geq 0}^n$ are given and that
$$ \sum_{i=1}^{n}{a_i \phi_j(x_i)} =  0 \qquad \mbox{for all}~ 0 \leq j \leq k.$$
We recall that the trivial eigenfunction $\phi_0$ is a constant and since it is exactly integrated, we obtain
$$ \sum_{i=1}^{n}{a_i} = \vol(M).$$
This guarantees that constants are integrated exactly. Moreover, we note that the heat equation preserves the integral for all $t>0$ and thus
$$ \int_{M}{  e^{t\Delta} \sum_{k=1}^{n}{ a_k \delta_{x_k}} dx} =  \vol(M).$$
In particular, subtracting the constant function 1 yields a function that has mean value 0. Moreover, using the fact that all but the trivial eigenfunction $\phi_0$ has mean value 0, we see for all
 $t \geq 0$
\begin{align*}
  \left(\sum_{i=1}^{n}{ e^{t\Delta} a_i \delta_{x_i}} \right) - 1 &= \sum_{\ell = 0}^{\infty}{ \left\langle   \left(\sum_{i=1}^{n}{ e^{t\Delta} a_i\delta_{x_i}} \right) - 1 , \phi_{\ell}\right\rangle \phi_{\ell}} \\
 &=\sum_{\ell = 1}^{\infty}{ \left\langle   \left(\sum_{i=1}^{n}{ e^{t\Delta} a_i\delta_{x_i}} \right) - 1 , \phi_{\ell}\right\rangle \phi_{\ell}} \\
 &=\sum_{\ell = 1}^{\infty}{ \left\langle   \left(\sum_{i=1}^{n}{ e^{t\Delta} a_i\delta_{x_i}} \right)  , \phi_{\ell}\right\rangle \phi_{\ell}}.
\end{align*}
We now use the self-adjointness of the heat propagator, the fact that $\phi_{\ell}$ are Laplacian eigenfunctions and the fact that these are evaluated exactly for the first $k$ eigenfunctions to conclude that
\begin{align*}
 \sum_{\ell = 1}^{\infty}{ \left\langle   \left(\sum_{i=1}^{n}{ e^{t\Delta} a_i\delta_{x_i}} \right)  , \phi_{\ell}\right\rangle \phi_{\ell}} &= 
 \sum_{\ell = 1}^{\infty}{ \left\langle  e^{t\Delta}  \left(\sum_{i=1}^{n}{ a_i\delta_{x_i}} \right)  , \phi_{\ell}\right\rangle \phi_{\ell}} \\
&= \sum_{\ell = 1}^{\infty}{ \left\langle   \left(\sum_{i=1}^{n}{  a_i\delta_{x_i}} \right)  ,e^{t\Delta} \phi_{\ell}\right\rangle \phi_{\ell}} \\
&= \sum_{\ell = 1}^{\infty}{ \left\langle   \left(\sum_{i=1}^{n}{ a_i\delta_{x_i}} \right)  ,e^{-\lambda_{\ell} t} \phi_{\ell}\right\rangle \phi_{\ell}}\\
&= \sum_{\ell = 1}^{\infty}{ e^{-\lambda_{\ell} t}  \left\langle \left(\sum_{i=1}^{n}{ a_i\delta_{x_i}} \right)  , \phi_{\ell}\right\rangle \phi_{\ell}}\\
&=  \sum_{\ell = k}^{\infty}{ \left\langle  \left(\sum_{i=1}^{n}{ a_i \delta_{x_i}} \right)  , e^{-\lambda_{\ell} t}  \phi_{\ell}\right\rangle \phi_{\ell}}.
\end{align*}
If we combine this with our observation from above, we obtain
$$ \frac{ \left\| \sum_{i=1}^{n}{ \nabla e^{t\Delta}a_i \delta_{x_i}}\right\|_{L^2}^2}{ \left\| \left(\sum_{i=1}^{n}{ e^{t\Delta}a_i \delta_{x_i}} \right) - 1\right\|_{L^2}^2} \geq \lambda_k.$$
The denominator is easy to analyze: observe that
$$
\left\| \left(\sum_{i=1}^{n}{ e^{t\Delta}a_i \delta_{x_i}} \right) - 1\right\|_{L^2}^2 = \int_{M}{  \left(\sum_{i=1}^{n}{ e^{t\Delta}a_i \delta_{x_i}} \right)^2 - 2  \left(\sum_{i=1}^{n}{ e^{t\Delta}a_i \delta_{x_i}} \right) + 1 ~dx} 
$$
Since the heat equation on compact manifolds preserves the mean value
$$ \int_{M}{ 2  \sum_{i=1}^{n}{ e^{t\Delta}a_i \delta_{x_i}}  dx} =  2\sum_{i=1}^{n}{a_i} = 2\vol(M)$$
and thus
$$ \left\| \left(\sum_{i=1}^{n}{ e^{t\Delta}a_i \delta_{x_i}} \right) - 1\right\|_{L^2}^2 = \left\| \sum_{i=1}^{n}{ e^{t\Delta}a_i \delta_{x_i}} \right\|_{L^2}^2 - \vol(M).$$
We will now bound the numerator from above
\begin{align*}
 \left\| \sum_{i=1}^{n}{ \nabla e^{t\Delta} a_i\delta_{x_i}}\right\|_{L^2}^2 = \sum_{i,j=1}^{n}{\left\langle \nabla e^{t \Delta} a_i \delta_{x_i}, \nabla e^{t\Delta}a_j\delta_{x_j}\right\rangle} 
\end{align*}
by invoking bounds on the logarithmic derivative of the heat kernel. These estimates were first obtained by Sheu \cite{sheu}, we also refer to Hsu \cite{hsu, hsu2} and Stroock \& Turetsky \cite{stroock}. These
bounds imply that, for $c>0$ depending only on $(M,g)$,
$$ \left| \frac{ \partial_y \left[ e^{t\Delta} \delta_x \right](y)}{ \left[ e^{t\Delta} \delta_x \right](y)} \right| = \left| \partial_y  \log{\left(\left[ e^{t\Delta} \delta_x \right](y)\right)} \right| \leq c\left(\frac{1}{\sqrt{t}} + \frac{d(x,y)}{t}\right).$$
Therefore $I:=  \left\langle \nabla e^{t \Delta} a_i \delta_{x_i}, \nabla e^{t\Delta}a_j\delta_{x_j}\right\rangle$ can be bounded by
\begin{align*}
I \leq c^2 \int_{M}{ \left(\frac{1}{\sqrt{t}} + \frac{d(x_i,x)}{t}\right) \left(\frac{1}{\sqrt{t}} + \frac{d(x_j,x)}{t}\right) \left[e^{t \Delta} a_i \delta_{x_i}\right](x)  \left[ e^{t\Delta}a_j\delta_{x_j}  \right](x) dx}.
\end{align*}
Suppose now that $d(x_i, x_j) \lesssim \sqrt{t}$. Then the Gaussian bounds on the kernel imply
$$ I \lesssim \frac{1}{t} \left\langle e^{t \Delta} \delta_{x_i}, e^{t\Delta} d_{x_j} \right\rangle,$$
where the implicit constant depends only on the manifold. Suppose now that $d(x_i, x_j) \gtrsim \sqrt{t}$. Then the dominant term is $d(x_i, x)/t \gtrsim t^{-1/2}$ and we may bound
$$I \lesssim  \int_{M}{  \frac{d(x_i,x)}{t} \frac{d(x_j,x)}{t} \left[e^{t \Delta} a_i \delta_{x_i}\right](x)  \left[ e^{t\Delta}a_j\delta_{x_j}  \right](x) dx}.$$
Gaussian bounds imply again that, for an implicit constant depending on the manifold
$$ I \lesssim   \frac{d(x_i,x_j)^2}{t^2} \int_{M}{ \left[e^{t \Delta} a_i \delta_{x_i}\right](x)  \left[ e^{t\Delta}a_j\delta_{x_j}  \right](x) dx}$$
which implies the result.
\end{proof}
A slightly more careful estimate implies shows that
$$ I \lesssim \frac{d(x_i, x_j)}{t^{3/2}} \int_{M}{ \left[e^{t \Delta} a_i \delta_{x_i}\right](x)  \left[ e^{t\Delta}a_j\delta_{x_j}  \right](x) dx}$$
but this has very little impact on the main result in the way that we use it since the contribution from these interactions are relatively small to begin with.
The proof also explain why the inequality should generally be very close to sharp as $t$ becomes large. The elementary fact that for functions of the type
$$ g(x) = \sum_{\ell =k}^{\infty}{ \left\langle g, \phi_{\ell}\right\rangle \phi_{\ell}} \quad \mbox{we have} \quad \left\| \nabla g\right\|_{L^2}^2 \geq \lambda_k \left\|g\right\|_{L^2}^2$$
is obviously sharp. Note, however, that our selection of $g = \sum_{i=1}^{n}{e^{t \Delta} \delta_{x_i}} - 1$ is far from arbitrary. Indeed, 
$$  \left\| \nabla \left(\sum_{i=1}^{n}{e^{t \Delta} \delta_{x_i}} -1\right) \right\|_{L^2}^2 = \left| \left\langle \Delta   \sum_{i =1}^{n}{ e^{t \Delta} \delta_{x_i}} ,     \sum_{i =1}^{n}{ e^{t \Delta} \delta_{x_i}} \right\rangle\right| $$
We now expand $\sum_{i=1}^{n}{ e^{t\Delta} \delta_{x_i}}$ into eigenfunctions and obtain
$$  \left| \left\langle \Delta   \sum_{i =1}^{n}{ e^{t \Delta} \delta_{x_i}} ,     \sum_{i =1}^{n}{ e^{t \Delta} \delta_{x_i}} \right\rangle\right| =  \sum_{\ell =k}^{\infty}{ \lambda_{\ell} e^{- 2\lambda_{\ell}t} \left|\left\langle  \sum_{i=1}^{n}{e^{t \Delta} \delta_{x_i}} , \phi_{\ell}\right\rangle\right|^2 }$$
while
$$  \left| \left\langle    \sum_{i =1}^{n}{ e^{t \Delta} \delta_{x_i}} ,     \sum_{i =1}^{n}{ e^{t \Delta} \delta_{x_i}} \right\rangle\right| =  \sum_{\ell =k}^{\infty}{ e^{- 2\lambda_{\ell}t} \left|\left\langle  \sum_{i=1}^{n}{e^{t \Delta} \delta_{x_i}} , \phi_{\ell}\right\rangle\right|^2 }$$
This implies that 
$$ \lim_{t \rightarrow \infty}{ \frac{ \left\| \nabla \left( \sum_{i=1}^{n}{e^{t \Delta} \delta_{x_i}} - 1 \right)  \right\|_{L^2}^2}{ \left\|  \sum_{i=1}^{n}{e^{t \Delta} \delta_{x_i}} - 1\right\|_{L^2}^2}} = \lambda_k.$$
It is also worth pointing out that the estimate on the logarithm of the heat kernel is on the full gradient in $(t,x,y)$ whereas we only use a special case
$$ \left|  \partial_y  \log{\left(\left[ e^{t\Delta} \delta_x \right](y)\right)} \right| \leq \left| \nabla_{x,y,t}  \log{\left(\left[ e^{t\Delta} \delta_x \right](y)\right)} \right|.$$
Moreover, the locally Euclidean case suggests that in our special case the term involving distance comes with a negative sign (this is a crucial ingredient in the proof of Theorem 2); furthermore, short-time asymptotics (we refer
again to the next section) show that this holds in the dominant short range regime. In any case, it certainly underlines that
 $$\mbox{simplified energy} = \sum_{i,j=1}^{n}{\frac{a_i a_j}{(8\pi t)^{d/2}} \exp\left( -\frac{d(x_i, x_j)^2}{8t}\right) },$$
 might be the quantity that controls most of the relevant structure and could be the most promising quantity for the purpose of numerical use.

\section{Proof of Theorem 2}
The key idea is to observe that in the Euclidean case, the explicit formula for the heat kernel can be used to show that
$$   \left\langle \nabla e^{t \Delta} \delta_{x_i}, \nabla e^{t\Delta}\delta_{x_j}\right\rangle   \leq \frac{\dim(M)}{4t} \left\langle e^{t \Delta} \delta_x, e^{t \Delta} \delta_y \right\rangle.$$
This is done as follows: we start by using that $\nabla, \Delta$ and $e^{t\Delta}$ are self-adjoint spectral multiplier that commute
$$  \left\langle \nabla e^{t \Delta} \delta_{x}, \nabla e^{t\Delta}\delta_{y}\right\rangle =   \left\langle - \Delta  e^{2t \Delta} \delta_{x},\delta_{y}\right\rangle.$$
Since the heat propagator, by definition, solves the heat equation,
$$ - \Delta  e^{2t \Delta} \delta_{x} = -\frac{1}{2} \partial_t e^{2t \Delta} \delta_{x}$$
and thus
\begin{align*}
-\frac{1}{2} \left[\partial_t e^{2t \Delta} \delta_{x}\right](y) &= - \frac{1}{2}\partial_t \frac{1}{(8 \pi t)^{d/2}} \exp\left( - \frac{d(x,y)^2}{8t} \right) \\
&=  \frac{1}{(8 \pi t)^{d/2} }   \exp\left( - \frac{d(x,y)^2}{8t} \right) \left( \frac{d}{4 t} - \frac{d(x,y)^2}{8t^2}\right) \\
&\leq \frac{d}{4t}  \frac{1}{(8 \pi t)^{d/2} }   \exp\left( - \frac{d(x,y)^2}{8t} \right) \\
&=  \frac{d}{4t} \left\langle e^{t \Delta} \delta_x,  e^{t \Delta} \delta_y \right\rangle.
\end{align*}
The main idea in the proof of Theorem 2 is to show that this inequality 'essentially' also works on general compact manifolds if $t$ is small. 
We first prove the desired result on $\mathbb{T}^2$ to outline the rather transparent proof and work with explicit constants.
We will then explain the modification necessary to make it work on a general $n-$dimensional manifold.

\begin{proof}[Proof on $\mathbb{T}^2$.] Note that
$$ \left[ e^{t\Delta} \delta_x\right](y) = \frac{1}{4\pi t} \sum_{k \in 2\pi\mathbb{Z}^2}{ e^{-\frac{\|x-y+k\|^2}{4t}}}.$$
Repeating the approach from above
$$   \left\langle \nabla e^{t \Delta} \delta_{x_i}, \nabla e^{t\Delta}\delta_{x_j}\right\rangle  =   \left\langle -\Delta e^{2 t \Delta} \delta_{x_i}, \delta_{x_j}\right\rangle = \left[- \frac{1}{2} \partial_t  e^{2 t \Delta} \delta_{x_i}\right](x_j).$$
We can now use the explicit form of the heat kernel to compute
\begin{align*}
\left[- \frac{1}{2} \partial_t  e^{2 t \Delta} \delta_{x_i}\right](x_j) &=  - \frac{1}{2} \partial_t  \frac{1}{8\pi t} \sum_{k \in 2\pi\mathbb{Z}^2}{ e^{-\frac{\|x_i-x_j+k\|^2}{8t}}} \\
&= \frac{1}{16 \pi t^2} \sum_{k \in 2\pi\mathbb{Z}^2}{ e^{-\frac{\|x_i-x_j+k\|^2}{8t}}} - \sum_{k \in 2\pi\mathbb{Z}^2}{\frac{\|x_i - x_j + k\|^2}{128 \pi t^3}  e^{-\frac{\|x_i-x_j+k\|^2}{8t}}}\\
&\leq  \frac{1}{2t}\frac{1}{8 \pi t} \sum_{k \in 2\pi\mathbb{Z}^2}{ e^{-\frac{\|x_i-x_j+k\|^2}{8t}}} = \frac{1}{2t} \left\langle e^{2t \Delta} \delta_{x_i}, \delta_{x_j} \right\rangle\\
&= \frac{1}{2t}  \left\langle e^{t \Delta} \delta_{x_i}, e^{t \Delta} \delta_{x_j} \right\rangle.
\end{align*}
The first step in the proof of the main result implies that
$$ \lambda_k \leq \frac{\sum_{i,j=1}^{n}{a_i a_j \left\langle \nabla e^{t\Delta} \delta_{x_i}, \nabla e^{t\Delta} \delta_{x_j} \right\rangle}}{ \sum_{i,j=1}^{n}{a_i a_j \left\langle  e^{t\Delta} \delta_{x_i},  e^{t\Delta} \delta_{x_j} \right\rangle} - \vol(M)}.$$
This can now be bounded from above by
\begin{align*}
 \lambda_k  &\leq  \frac{1}{2t} \frac{\sum_{i,j=1}^{n}{a_i a_j \left\langle  e^{t\Delta} \delta_{x_i},  e^{t\Delta} \delta_{x_j} \right\rangle}}{ \sum_{i,j=1}^{n}{a_i a_j \left\langle  e^{t\Delta} \delta_{x_i},  e^{t\Delta} \delta_{x_j} \right\rangle} - 4\pi^2}\\
&\leq \frac{1}{2t} \left[ 1 + \frac{4 \pi^2}{ \sum_{i,j=1}^{n}{a_i a_j \left\langle  e^{t\Delta} \delta_{x_i},  e^{t\Delta} \delta_{x_j} \right\rangle} - 4\pi^2} \right]
\end{align*}
We conclude by observing that since all the weights are nonnegative
$$  \sum_{i,j=1}^{n}{a_i a_j \left\langle  e^{t\Delta} \delta_{x_i},  e^{t\Delta} \delta_{x_j} \right\rangle} \geq   \sum_{i=1}^{n}{a_i^2 \left\langle  e^{t\Delta} \delta_{x_i},  e^{t\Delta} \delta_{x_i} \right\rangle}.$$
Clearly,
\begin{align*}
 \left\langle  e^{t\Delta} \delta_{x_i},  e^{t\Delta} \delta_{x_i} \right\rangle =   \left\langle  e^{2t\Delta} \delta_{x_i},  \delta_{x_i} \right\rangle = \frac{1}{8\pi t} \sum_{k \in 2\pi\mathbb{Z}^2}{ e^{-\frac{\|k\|^2}{8t}}} \geq \frac{1}{8 \pi t}.
\end{align*}
Since the quadrature formula integrates constants exactly,
$$ 4 \pi^2 = \sum_{i=1}^{n}{a_i} \leq \sqrt{n} \left( \sum_{i=1}^{n}{a_i^2} \right)^{1/2} \quad \mbox{and thus} \quad    \sum_{i=1}^{n}{a_i^2}  \geq \frac{16 \pi^4}{n}.$$
Altogether, this implies
 $$   \sum_{i,j=1}^{n}{a_i a_j \left\langle  e^{t\Delta} \delta_{x_i},  e^{t\Delta} \delta_{x_j} \right\rangle} \geq \frac{1}{8 \pi t}  \frac{16 \pi^4}{n}$$
and thus
$$ \lambda_k \leq   \frac{1}{2t} \frac{ \frac{1}{8 \pi t}  \frac{16 \pi^4}{n}}{\frac{1}{8 \pi t}  \frac{16 \pi^4}{n} - 4 \pi^2}.$$
Optimization in $t$ suggests to pick $t = \pi/(4n)$ which then yields
$$ \lambda_k \leq \frac{4n}{\pi}.$$
The Weyl law implies that on $\mathbb{T}^2 \cong [0,2\pi]^2$
$$ \lambda_k = \frac{k}{\pi} + o(k)$$
and thus, for $k$ sufficiently large,
$$ k \leq 4n + o(n).$$
\end{proof}

\begin{proof}[Proof of general manifolds]
We can assume w.l.o.g. $\vol(M) = 1$. We proceed as before and note that
$$ 1 = \left(\sum_{i=1}^{n}{a_i}\right)^2 \leq  n \sum_{i=1}^{n}{a_i^2}$$
and therefore
$$ \sum_{i=1}^{n}{a_i^2} \geq \frac{1}{n}.$$
It remains to estimate the interaction of the heat kernel, in particular we require lower bounds on
$ \left\langle e^{t\Delta} \delta_x, e^{t\Delta} \delta_x\right\rangle$ and upper bounds on
 $ \left\langle \nabla e^{t\Delta} \delta_x, \nabla e^{t\Delta} \delta_y\right\rangle$. As before,
we can rewrite the second expression as
$$  \left\langle \nabla e^{t\Delta} \delta_x, \nabla e^{t\Delta} \delta_y\right\rangle = \left[-\Delta e^{2 t \Delta} \delta_{x_i}\right](x_j) = \left[- \frac{1}{2} \partial_t  e^{2 t \Delta} \delta_{x_i}\right](x_j).$$
We use the asymptotic expansion (see e.g. \cite[Theorem 5.1.1.]{hsu2})
$$ \left[ e^{t \Delta}\delta_x\right](y) \sim  \left(\frac{1}{4\pi t}\right)^{\frac{d}{2}} e^{-\frac{d(x,y)^2}{4t}} \sum_{n=0}^{\infty}{H_n(x,y) t^n}$$
which is valid uniformly as $t \rightarrow 0$ on any compact subset away from the cut locus. $H_n$ are smooth functions, $H_0(x,y) > 0$ and $H_0(x,x) = 1$. As before, 
\begin{align*}
  \sum_{i,j=1}^{n}{a_i a_j \left\langle e^{t\Delta} \delta_{x_i}, e^{t \Delta} \delta_{x_j} \right\rangle } &\geq  \sum_{i=1}^{n}{a_i a_j \left\langle e^{t\Delta} \delta_{x_i}, e^{t \Delta} \delta_{x_j} \right\rangle }\\
&\geq \left( \inf_{y \in M} \left\langle e^{2t \Delta} \delta_y, \delta_y \right\rangle \right) \sum_{i=1}^{n}{a_i^2} \\
&\geq \left(1 + \mathcal{O}(t)\right)\frac{1}{(8 \pi t)^{d/2}} \frac{1}{n}.
\end{align*}
We will work on time scale $t \sim n^{-2/d}$. On that scale, the dominant local interactions are local and at scale $\sim n^{-1/d}$ and we are able to invoke short-time asymptotics.
 In particular
$$ -\frac{1}{2} \partial_t \left[ e^{2t \Delta}\delta_x\right](y) = \partial_t  \left[ \left(\frac{1}{8\pi t}\right)^{\frac{d}{2}} e^{-\frac{d(x,y)^2}{8t}} \right]  \sum_{n=0}^{\infty}{H_n(x,y)(2t)^n} + \mbox{error},$$
where the derivative acting on the bracket reproduces exactly the Euclidean behavior and yields a quantity that is $\sim t^{-1}$ times as large as the pure heat kernel interaction. The error is at the
order of the heat kernel interaction itself
$$ \left|\mbox{error}\right|  \lesssim \left(\frac{1}{8\pi t}\right)^{\frac{d}{2}} e^{-\frac{d(x,y)^2}{8t}}.$$
The next step consists of putting interactions $ \left\langle \nabla e^{t\Delta} \delta_x, \nabla e^{t\Delta} \delta_y\right\rangle$
into two different regimes: those for which $d(x_i, x_j) \leq n^{-1/10d}$ and those for which $d(x_i, x_j) \geq n^{-1/10d}$. Any pair in the first group is, for $n$ sufficiently
large, in the regime where the asymptotic expansion applies (because the manifold is compact and thus has an injectivity radius uniformly bounded away from 0). For the second pair, it suffices to apply the
bounds on the logarithmic derivative to conclude that
\begin{align*}
 \left| \left\langle \nabla e^{t\Delta} \delta_x, \nabla e^{t\Delta} \delta_y\right\rangle \right| &\lesssim c\frac{d(x,y)^2}{t^2}  \left\langle e^{t\Delta} \delta_x, e^{t\Delta} \delta_y\right\rangle.
\end{align*}
However, using standard Gaussian bounds on the heat kernel, 
\begin{align*}
 \sum_{d(x_i, x_j) \geq n^{-\frac{1}{10d}}}{ \frac{d(x,y)^2}{t^2}  \left\langle e^{t\Delta} \delta_x, e^{t\Delta} \delta_y\right\rangle} &\lesssim \frac{1}{t^2} \sum_{d(x_i, x_j) \geq n^{-\frac{1}{10d}}}{ \left\langle e^{t\Delta} \delta_x, e^{t\Delta} \delta_y\right\rangle} \\
&\lesssim \frac{1}{t^2} \frac{1}{t^{d/2}} \sum_{d(x_i, x_j) \geq n^{-\frac{1}{10d}}}{ \exp\left(-c \frac{d(x,y)^2}{t} \right)} \\
&\lesssim  \frac{n^2}{t^2} \frac{1}{t^{d/2}}  \exp\left(-\frac{c}{n^{\frac{1}{5d}} t} \right)
\end{align*}
Since $t \sim n^{-2/d}$, this term is superpolynomially decreasing in $n$ and thus
$$  \sum_{d(x_i, x_j) \geq n^{-\frac{1}{10d}}}{ \left| \left\langle \nabla e^{t\Delta} \delta_x, \nabla e^{t\Delta} \delta_y\right\rangle \right|} \ll \sum_{i,j=1}^{n}{  \left| \left\langle  e^{t\Delta} \delta_x,  e^{t\Delta} \delta_y\right\rangle \right|}.$$
As a consequence, we obtain
\begin{align*}
\sum_{i,j=1}^{n}{a_i a_j \left\langle \nabla e^{t\Delta} \delta_{x_i}, \nabla e^{t\Delta} \delta_{x_j} \right\rangle} &\leq
\frac{d}{4t}\sum_{i,j=1}^{n}{a_i a_j \left\langle e^{t\Delta} \delta_{x_i}, e^{t\Delta} \delta_{x_j} \right\rangle} \\
&+ c\sum_{i,j=1}^{n}{a_i a_j \left\langle  e^{t\Delta} \delta_{x_i},  e^{t\Delta} \delta_{x_j} \right\rangle},
\end{align*}
where $c$ only depends on the manifold.
Then, however,
$$ \lambda_k \leq   \frac{d}{4t} \frac{ \left(1 + (4 c/d) t\right) \sum_{i,j=1}^{n}{a_i a_j \left\langle e^{t\Delta} \delta_{x_i}, e^{t\Delta} \delta_{x_j} \right\rangle} }{\sum_{i,j=1}^{n}{a_i a_j \left\langle e^{t\Delta} \delta_{x_i}, e^{t\Delta} \delta_{x_j} \right\rangle} - 1},$$
which is monotonically decreasing in the sum term. As outlined above, we can bound the sum from below by 
\begin{align*}
  \sum_{i,j=1}^{n}{a_i a_j \left\langle e^{t\Delta} \delta_{x_i}, e^{t\Delta} \delta_{x_j} \right\rangle} &\geq  \sum_{i=1}^{n}{a_i^2 \left\langle e^{t\Delta} \delta_{x_i}, e^{t\Delta} \delta_{x_i} \right\rangle} \\
&\geq  \left(1 + \mathcal{O}(t)\right)\frac{1}{(8 \pi t)^{d/2}} \frac{1}{n}.
\end{align*}
 and thus, for $t \sim  n^{-2/d}$
$$ \lambda_k \leq   \frac{d}{4t} \frac{ \left(1 +  (4c/d)  t\right) (8 \pi t)^{-d/2} n^{-1}}{  (8 \pi t)^{-d/2} n^{-1} - 1} \sim \frac{1}{t} t^{-d/2} n^{-1} \sim \frac{1}{t} \sim n^{2/d}.$$
Since $\vol(M)=1$, Weyl's law implies that $ \lambda_k \sim k^{2/d}$ and thus $k \lesssim n$. 
The computation of the constant $ c_d$ is a simple consequence of being slightly more careful: by minimizing the bound
$$ \lambda_k \leq   \frac{d}{4t} \frac{ \left(1 + c_2 c_1^{-1} t\right) (8 \pi t)^{-d/2} n^{-1}}{  (8 \pi t)^{-d/2} n^{-1} - 1}$$
in $t$, we obtain
$$ \lambda_k \leq 2^{1 - \frac{2}{d}} (d+2)^{\frac{2}{d}+1} \pi n^{2/d} + o(n^{2/d}).$$
(For $d=2$ this differs from the result on $\mathbb{T}^2$ above because we now work with the normalization $\vol(M) = 1$).
The sharp constant in Weyl's law on a manifold with $\vol(M)=1$ is
$$ \lambda_k^{d/2} \sim \frac{(2\pi)^d}{\omega_d} k^{}  + o(k)  = \frac{     (2\pi)^d  }{\pi^{\frac{d}{2}}}  \Gamma\left( \frac{d}{2} + 1\right)  k^{} + o(k),$$
where $\omega_d$ is, as usual, the volume of the unit ball in $\mathbb{R}^d$. Combined, this yields 
$$ k \leq  \frac{\left( \frac{d}{2} + 1\right)^{\frac{d}{2} + 1}}{\Gamma\left(\frac{d}{2} + 1\right)}  n + o(n).$$
\end{proof}

There are two obvious spots where the argument could be improved. The first one is in the use of the lossy estimate
$$  \sum_{i,j=1}^{n}{a_i a_j \left\langle  e^{t\Delta} \delta_{x_i},  e^{t\Delta} \delta_{x_j} \right\rangle} \geq   \sum_{i=1}^{n}{a_i^2 \left\langle  e^{t\Delta} \delta_{x_i},  e^{t\Delta} \delta_{x_i} \right\rangle},$$
which completely ignores off-diagonal contributions. Since we are already working at spatial scale $\sim n^{-1/d}$, the actual improvement in the constant is likely to be small. A much more substantial improvement,
especially in higher dimensions, is likely to follow from operating on larger time scales (which will of course require an understanding of local interactions). The main reason for assuming
this to be the case has already been given after the proof of Theorem 1: the heat equation is actively suppressing higher frequencies at a faster rate and this effect becomes more pronounced as time
becomes large (at the cost of a greater combinatorial complexity).\\

\textbf{Acknowledgment.} The author is grateful to Vladimir Rokhlin for several insightful conversations about quadratures.

\end{document}